\documentclass[12pt,a4paper]{article}
\usepackage{indentfirst}
\usepackage{amsmath,amssymb,amsfonts,amsthm,graphics}
\usepackage{makeidx}
\newtheorem{thm}{Theorem}[section]

\theoremstyle{definition}

\def\-{\mbox{--}}

\begin{document}

\title{\bf \Large Note on rainbow connection number \\
of dense graphs\footnote{Supported by NSFC No.11071130.} }
\author{\small Jiuying Dong, Xueliang Li\\
\small Center for Combinatorics and LPMC-TJKLC \\
\small Nankai University, Tianjin 300071, China\\
\small Email: jiuyingdong@126.com; lxl@nankai.edu.cn}

\date{}

\maketitle{}

\begin{abstract}

An edge-colored graph $G$ is rainbow connected if any two vertices
are connected by a path whose edges have distinct colors. The
rainbow connection number of a connected graph $G$, denoted by
$rc(G)$, is the smallest number of colors that are needed in order
to make $G$ rainbow connected. Following an idea of Caro et al., in
this paper we also investigate the rainbow connection number of
dense graphs. We show that for $k\geq 2$, if $G$ is a non-complete
graph of order $n$ with minimum degree $\delta (G)\geq
\frac{n}{2}-1+log_{k}{n}$, or minimum degree-sum $\sigma_{2}(G)\geq
n-2+2log_{k}{n}$, then $rc(G)\leq k$; if $G$ is a graph of order $n$
with diameter 2 and $\delta (G)\geq
2(1+log_{\frac{k^{2}}{3k-2}}{k})log_{k}{n}$, then $rc(G)\leq k$. We
also show that if $G$ is a non-complete bipartite graph of order $n$
and any two vertices in the same vertex class have at least
$2log_{\frac{k^{2}}{3k-2}}{k}log_{k}{n}$ common neighbors in the
other class, then $rc(G)\leq k$.\\[3mm]
\end{abstract}

\section{ Introduction}

All graphs under our consideration are finite, undirected and
simple. For notation and terminology not defined here, we refer to
\cite{bondy2008graph}. Let $G$ be a graph. The length of a path in
$G$ is the number of edges of the path. The distance between two
vertices $u$ and $v$ in $G$, denoted by $d(u, v)$, is the length of
a shortest path connecting them in $G$. If there is no path
connecting $u$ and $v$, we set $d(x,y):=\infty$. An edge-coloring of
a graph is a function from its edges set to the set of natural
numbers. A graph $G$ is rainbow edge-connected if for every pair of
distinct vertices $u$ and $v$ of $G$, $G$ has a $u-v$ path whose
edges are colored with distinct colors. This concept was introduced
by Chartrand et al. \cite{chartrand2008rainbow}. The minimum number
of colors required to rainbow color a connected graph is called its
rainbow connection number, denoted by $rc(G)$. Observe that if $G$
has $n$ vertices, then $rc(G)\leq n-1$. Clearly, $rc(G) \geq
diam(G)$, the diameter of $G$. In \cite{chartrand2008rainbow},
Chartrand et al. determined the rainbow connection numbers of
wheels, complete graphs and all complete multipartite graphs. In
\cite {chakraborty2009hardness}, Chakraborty et al. proved that
given a graph $G$, deciding if $rc(G)=2$ is NP-Complete. In
particular, computing $rc(G)$ is NP-Hard.

If $\delta (G)\geq \frac{n}{2}$, then $diam(G)=2$, but we do not
know if this guarantees $rc(G)=2$. In \cite{caro2008rainbow}, Caro
et al. investigated the rainbow connection number of dense graphs,
and they got the following results.
\begin{thm}
Any non-complete graph with $\delta (G)\geq \frac{n}{2}+logn$ has
$rc(G)=2$.
\end{thm}
\begin{thm}
Let $c=\frac{1}{log(9/7)}$. If $G$ is a non-complete bipartite graph
with $n$ vertices and any two vertices in the same vertex class have
at least $2clogn$ common neighbors in the other class, then
$rc(G)=3$.
\end{thm}

We will follow their idea to investigate dense graphs again. And we
get the following results.
\begin{thm}
Let $k\geq 2$ be an integer. If $G$ is a non-complete graph of order
$n$ with $\delta (G)\geq \frac{n}{2}-1+log_{k}{n}$, then $rc(G)\leq
k$.
\end{thm}
\begin{thm}
Let $k\geq 2$ be an integer. If $G$ is a non-complete graph of order
$n$ with $\sigma_{2}(G)\geq n-2+2log_{k}{n}$, then $rc(G)\leq k$.
\end{thm}
\begin{thm}
Let $k\geq 3$ be an integer. If $G$ is a non-complete bipartite
graph of order $n$ and any two vertices in the same vertex class
have at least $2log_{\frac{k^{2}}{3k-2}}klog_{k}n$ common
neighbors in the other class, then $rc(G)\leq k$.
\end{thm}

In \cite {chakraborty2009hardness}, Chakraborty et al. showed the
following result.
\begin{thm}
If $G$ is a graph of order $n$ with diameter 2 and $\delta (G)\geq
8logn$, then $rc(G)\leq 3$. Furthermore, such a coloring is given
with high probability by a uniformly random 3-edge-coloring of the
graph $G$, and can also be found by a polynomial time deterministic
algorithm.
\end{thm}

Now we get the following result.
\begin{thm}
Let $k\geq 3$ be an integer. If $G$ is a graph of order $n$ with
diameter 2 and $\delta (G)\geq
2(1+log_{\frac{k^{2}}{3k-2}}k)log_{k}n$, then $rc(G)\leq k$.
\end{thm}

\section{ Proof of the theorems }

{\bf Proof of Theorem 1.3:} Let $G$ be a non-complete  graph of
order $n$ with $\delta (G)\geq \frac{n}{2}-1+log_{k}n$. We use
$k$ different colors to randomly color every edge of $G$. In the
following we will show that with positive probability, such a random
coloring make $G$ rainbow connected. For any pair $u,v\in V(G),
uv\not\in E(G)$, since $d(u)\geq \frac{n}{2}-1+log_{k}n$,
$d(v)\geq \frac{n}{2}-1+log_{k}n$, there are at least
$2log_{k}n$ common neighbors between $u$ and $v$, that is
$|N(u)\cap N(v)|\geq 2log_{k}n$. Hence  there are at least
$2log_{k}n$ edge-disjoint paths of length two from $u$ to $v$.
For any $w\in N(u)\cap N(v)$, the probability that the path $uwv$ is
not a rainbow path is $\frac{1}{k}$. Hence, the probability that all
these edge-disjoint paths are not rainbow is at most
$(\frac{1}{k})^{2log_{k}n}=\frac{1}{n^{2}}$. Since there are less
than $n\choose 2$ pairs non-adjacent vertices, and $ {n\choose
2}\frac{1}{n^{2}} < 1$. We may get that with positive probability,
each pair of non-adjacent vertices are connected by a rainbow path.
This completes the proof of Theorem 1.3.

{\bf Proof of Theorem 1.4:} Let $G$ be a non-complete  graph of
order $n$ with $\sigma_{2}(G)\geq n-2+2log_{k}n$. We use $k$
different colors to randomly color every edge of $G$. In the
following we will show that with positive probability, such a random
coloring make $G$ rainbow connected. For any pair $u,v\in V(G),
uv\not\in E(G)$, as $\sigma_{2}(G)\geq n-2+2log_{k}n$,  it
follows that $|N(u)\cap N(v)|\geq 2log_{k}n$. Similar to  the
proof of Theorem 1.3, we may get that with positive probability,
each pair of non-adjacent vertices are connected by a rainbow path.
This completes the proof of Theorem 1.4.

{\bf Proof of Theorem 1.5:} Let $G$ be a non-complete bipartite
graph of order $n$ and any two vertices in the same vertex class
have at least $2log_{\frac{k^{2}}{3k-2}}klog_{k}n$ common
neighbors in the other class. We use $k$ different colors to
randomly color every edge of $G$.  In the following we will show
that with positive probability, such a random coloring make $G$
rainbow connected.  For every pair $u,v\in V(G)$ and  $u,v$ are in
the same class of $V(G)$, then the distance of $d(u,v)=2$, as
$|N(u)\cap N(v)|\geq 2log_{\frac{k^{2}}{3k-2}}klog_{k}n$,
there are at least $2log_{\frac{k^{2}}{3k-2}}klog_{k}n$
edge-disjoint paths of length two from $u$ to $v$. The probability
that all these edge-disjoint paths are not rainbow is at most
$(\frac{1}{k})^{2log_{\frac{k^{2}}{3k-2}}klog_{k}n}<(\frac{1}{k})^{2log_{k}n}=\frac{1}{n^{2}}$.
For every pair $u,v\in V(G)$ from different classes of $G$ and
$uv\not\in E(G)$, then the distance of $d(u,v)$ is 3. Fix a neighbor
$w_u$ of $u$, for any $u_i\in N(w_u)\cap N(v)$, the probability that
$uw_uu_iv$ is not a rainbow path is $\frac{3k-2}{k^{2}}$. We know
$|N(w_u)\cap N(v)|\geq 2log_{\frac{k^{2}}{3k-2}}klog_{k}n$.
Hence, the probability that all these edge-disjoint paths are not
rainbow is at most
$(\frac{3k-2}{k^{2}})^{2log_{\frac{k^{2}}{3k-2}}klog_{k}n}=\frac{1}{n^{2}}$.
Thus, we may get that with positive probability, each pair of
non-adjacent vertices are connected by a rainbow path. This
completes the proof of Theorem 1.5.

{\bf Proof of Theorem 1.7:} \\
Let $G$ be a graph of order $n$ with diameter 2.  We use $k$
different colors to randomly color every edge of $G$. In the
following we will show that with positive probability, such a random
coloring make $G$ rainbow connected. For any two non-adjacent
vertices $u,v$, if $|N(u)\cap N(v)|\geq 2log_{k}n$, then there
are at least $2log_{k}n$ edge-disjoint paths of length two from
$u$ to $v$. The probability that all these edge-disjoint paths are
not rainbow is at most
$(\frac{1}{k})^{2log_{k}n}=\frac{1}{n^{2}}$. Otherwise,
$|N(u)\cap N(v)|< 2log_{k}n$. Let $A=N(u)\setminus N(v),
B=N(v)\setminus N(u)$, then $|A|\geq
2log_{\frac{k^{2}}{3k-2}}klog_{k}n, |B|\geq
2log_{\frac{k^{2}}{3k-2}}klog_{k}n$.  As the diameter of $G$
is two, for any $x\in A$, $\exists y_{x}\in N(v)$ such that $xy_x\in
E(G)$, that is $xy_xv$ is a path of length 2.  Now, we will consider
the set of at least $ 2log_{\frac{k^{2}}{3k-2}}klog_{k}n$
edge-disjoint paths $P=\{uxy_xv: x\in A\}$.  For every $x\in A$, the
probability that $uxy_xv$ is not a rainbow path is
$\frac{3k-2}{k^{2}}$. Moreover, this
event is independent of the corresponding events for all other
members of $A$, because this probability does not change even with
full knowledge of the colors of all edges incident with $v$.
Therefore, the probability that all these edge-disjoint paths are
not rainbow is at most
$(\frac{3k-2}{k^{2}})^{2log_{\frac{k^{2}}{3k-2}}klog_{k}n}=\frac{1}{n^{2}}$.
Since there are less than $n\choose 2$  pairs non-adjacent vertices,
and $ {n\choose 2}\frac{1}{n^{2}} < 1$. We may get that with
positive probability, each pair of non-adjacent vertices are
connected by a rainbow path. This completes the proof of Theorem
1.7.

\end{document}